\documentclass[12pt,reqno]{article}
mm \textheight=8.9in

\usepackage{amssymb}
\usepackage{amsmath}
\usepackage{amsthm}

\usepackage{color}
\newcommand{\tw}[3]{{$#1$}${\,\scriptscriptstyle {#2}}\atop\raise9pt\hbox{$\scriptstyle\tp$} ${$#3$}}
\newcommand{\st}[1]{\mbox{${\,\scriptscriptstyle {#1}}\atop\raise5.5pt\hbox{$*$}$}}
\newcommand{\btr}{\raise1.2pt\hbox{$\scriptstyle\blacktriangleright$}\hspace{2pt}}
\newcommand{\id}{\mathrm{id}}

\newcommand{\tp}{\otimes}

\newcommand{\D}{\mathfrak{D}}

\newcommand{\Z}{\mathbb{Z}}

\newcommand{\la}{\lambda}

\newcommand{\End}{\mathrm{End}}

\newcommand{\Nc}{\mathcal{N}}

\newcommand{\tr}{\triangleright}

\newcommand{\btl}{\mbox{\raise1.1pt\hbox{$\scriptstyle\blacktriangleleft$}}}
\newcommand{\ad}{\mathrm{ad}}

\newcommand{\g}{\mathfrak{g}}

\newcommand{\h}{\mathfrak{h}}

\newcommand{\cth}{\mathrm{coth}}
\newcommand{\eps}{\epsilon}
\newcommand{\Alt}{\mathrm{Alt}}
\newcommand{\CYB}{\mathrm{CYB}}
\newcommand{\Exp}{\mathrm{Exp}\hspace{0.5pt}}
\newcommand{\nn}{\nonumber}

\renewcommand{\l}{\mathfrak{l}}

\newcommand{\al}{\alpha}

\newcommand{\be}{\begin{eqnarray}}
\newcommand{\ee}{\end{eqnarray}}

\newtheorem{thm}{Theorem}[section]
\newtheorem{propn}[thm]{Proposition}

\newtheorem{corollary}[thm]{Corollary}

\theoremstyle{definition}
\newtheorem{remark}[thm]{Remark}
\newtheorem{definition}[thm]{Definition}

\begin{document}
\title{Trigonometric dynamical r-matrices over Poisson Lie base\footnote{
This research is supported in part
by the Israel Academy of Sciences grant no. 8007/99-03,
the Emmy Noether Research Institute for Mathematics,
the Minerva Foundation of Germany,  the Excellency Center "Group
Theoretic Methods in the study of Algebraic Varieties"  of the Israel
Science foundation, and by the RFBR grant no. 03-01-00593.
} }
\author{A. Mudrov}
\date{}
\maketitle
\begin{center}
{Department of Mathematics, Bar Ilan University,  52900 Ramat Gan,
Israel,\\
Max-Planck Institut f$\ddot{\rm u}$r Mathematik, Vivatsgasse 7, D-53111 Bonn, Germany.\\
e-mail: mudrov@mpim-bonn.mpg.de}
\end{center}
\begin{abstract}
Let $\g$ be a finite dimensional complex Lie algebra and $\l\subset \g$ a
Lie subalgebra equipped with the structure of a factorizable quasitriangular Lie bialgebra.
Consider the Lie group $\Exp \l$ with the Semenov-Tjan-Shansky Poisson bracket
as a Poisson Lie manifold for the double Lie bialgebra $\D\l$.
Let $\Nc_\l(0)\subset \l$ be an open domain parameterizing a neighborhood of
the identity in $\Exp \l$ by the exponential map.
We present dynamical $r$-matrices with values in $\g\wedge \g$ over the Poisson Lie base manifold $\Nc_\l(0)$.
\end{abstract}
\maketitle
\section {Introduction}
The classical dynamical Yang-Baxter equation (CDYBE) is a differential equation which
generalizes the algebraic or ordinary classical Yang-Baxter equation (CYBE) associated with an arbitrary
Lie algebra $\g$. Contrary to the CYBE, a solution of CDYBE is a function of the so called
dynamical parameter. Up to recently, the space of parameters was taken to be
the dual $\l^*$ of a Lie subalgebra $\l$ in $\g$. In \cite{DM1,DM2}, the
DYBE was formulated for any quadruple $(\g,\l,\l^*,L)$, where
$\l\subset \g$ is a pair of Lie algebras, $(\l,\l^*)$ is a Lie bialgebra, and
$L$ is a special Poisson Lie (PL) manifold, the space of parameters.

A particular case of the PL version of the CDYBE for a quasitriangular Lie bialgebra $\l=\g$
and $L$ being the group space $\Exp \g$ appeared in \cite{BFP} in connection
with a factorization problem in the chiral WZW model. The corresponding dynamical
r-matrix was written out in \cite{FM}. In the present paper
we give examples of dynamical r-matrices for the quadruple $(\g,\l,\l^*,L)$, where
$\g$ is a complex Lie algebra and $(\l,\l^*)$ a factorizable quasitriangular
Lie bialgebra. As a base manifold $L$ we take a domain $\Nc_\l(0)$ in $\l$
parameterizing a neighborhood of the identity in the group $\Exp \l$
equipped with the Semenov-Tjan-Shansky (STS) Poisson bracket, \cite{S}.

We apply the Etingof-Varchenko approach of  base reduction, \cite{EV},
to the PL CDYBE.
Then we employ, within the PL setting, the idea of Etingof-Schiffmann which is
 used for constructing generalized Alekseev-Meinrenken dynamical r-matrices in \cite{ES}.
 In this way we obtain PL dynamical r-matrices that are in between the r-matrices
 of \cite{FM} and trigonometric r-matrices of \cite{EV,S}.

\section {Dynamical Yang-Baxter equation on PL base manifold}
Let $(\l,\l^*)$ be a finite dimensional Lie bialgebra and let $\delta$ denote the cobracket $\l\to \l\wedge \l$, \cite{D}.
Let $\D\l=\l\bowtie \l^*_{op}$ be the double
Lie algebra with the canonical invariant symmetric tensor $\theta:=\sum_i (\eta^i\tp \xi_i+\xi_i\tp \eta^i)$,
where $\{\xi_i\}\subset \l$ and $\{\eta^i\}\subset \l^*$ are the dual bases.
Recall from \cite{DM1} that a $\D\l$-manifold is called an $\l$-base manifold if the Casimir element
$\theta$ generates the zero bidifferential operator via the action of $\D\l$. It follows that
a base manifold is equipped with a $\D\h$-PL structure induced by the r-matrix
$\sum_i \eta^i\tp \xi_i\in (\D\l)^{\tp 2}$ via the action of $\D\l$.

In the present paper,  by a function on $L$ we understand an analytical or meromorphic function.
Let $\g$ be a Lie algebra containing $\l$ as a subalgebra.
Let us call a function $r\colon L\to \wedge^2\g$ quasi-invariant if
\be
 \xi\tr   r(\la) + [\xi\tp 1+ 1\tp \xi,  r(\la)] &=& -\delta(\xi), \quad \xi\in\l,
\label{inv}
\ee
where $\delta$ is the cobracket on $\l$ and $\tr$ denotes the $\l$-action on functions on $L$ by vector fields.
If the  Lie  bialgebra $\l$ is coboundary, i.e.
$\delta(\xi):=[\xi\tp 1+1\tp\xi, r_\l]$ for some $r_\l\in \wedge^2\l$,
then a function $r\colon L\to \wedge^2\g$ is quasi-invariant if and only if
the function $r':=r+r_\l$ is invariant, i.e.
\be
 \xi\tr   r'(\la) + [\xi\tp 1+ 1\tp \xi,  r'(\la)] &=&0 , \quad \xi\in\l.
\label{inv1}
\ee
\begin{definition}
\label{def_drm}
A quasi-invariant function $r\colon L\to \wedge^2\g$ is called a dynamical r-matrix over the base $L$ if
\be
\label{cDYBE}
\sum_{i} \mathrm{Alt}\bigl(\xi_i\tp  \nabla_{\eta^i}
r(\la)\bigr)+\mathrm{CYB}\bigl(r(\la)\bigr)
 &=& Z ,
\ee
where $\nabla_{\eta}$ is the vector field on $L$ generated by $\eta\in \l^*_{op}$, $Z\subset (\wedge^3 \g)^\g$
is  an invariant element, and
\be
\CYB(A)&:=& [A_{12},A_{13}]+[A_{12},A_{23}]+[A_{13},A_{23}], \quad A \in \g\tp \g,
\nn
\\
\Alt(B)&:=&B_{123}-B_{213}+B_{231}, \quad B \in \g\tp \g\tp \g.
\nn
\ee
\end{definition}
Note that Definition \ref{def_drm} differs from the definition given in \cite{DM1,DM2}
by the change of sign $r\to -r$.
Sometimes we will use the term PL dynamical r-matrix in order to distinguish it
from the particular case of abelian $\l^*$.

\begin{remark}
\label{dilation}
Suppose  $(\l,\l^*)$ is the trivial bialgebra, with zero $\delta$ (abelian $\l^*$).
Then replacement $Z$ by $\eps^2 Z$ results in the transformation
$r(\la)\to \eps r(\eps\la)$ of the dynamical r-matrix, the solution to (\ref{cDYBE}).
\end{remark}

Equation (\ref{cDYBE}) appeared in \cite{DM1,DM2} in connection with equivariant deformation quantization on
$G$-spaces. Its particular case  for $\l=\g$ and $L=\Exp \l$ was considered in
an equivalent form in the paper \cite{BFP} devoted to the chiral WZW model.

\section {Poisson Lie CDYBE on group manifold}
\label{GB}
Suppose that the Lie bialgebra $\l$ is quasitriangular, i.e.
coboundary with an r-matrix $r_\l\in \wedge^2 \l$ and an element $\Omega_\l\in (S^2\l)^\l$ such that
$\CYB(r_\l+\frac{1}{2}\Omega_\l)=0$.
Consider the connected simply connected Lie group $\Exp \l$ corresponding to the Lie bialgebra $\l$
and put $L=\Exp \l$.
According to \cite{S}, there is a PL structure on $L$ considered as an $\Exp \l$ -manifold
with respect to conjugation. In fact, $L$-is an $\l$-base manifold, \cite{DM1}.
The action of double Lie algebra $\D\l$ on $L$ is defined as follows.
For $\xi\in \l$ denote by $\xi^l$ and $\xi^r$ the vector fields
acting by
$$
\xi^l f(g)=\frac{d}{dt}f(ge^{t\xi})|_0
,\quad
\xi^r f(g)=\frac{d}{dt}f(e^{t\xi}g)|_0,\quad f\in C^\infty(L).
$$
Then $\D\h$ acts on $L$  by
the vector fields
\be
\xi^l-\xi^r \quad\mbox{for}\quad
\xi\in\l
 \quad\mbox{and}\quad
r^r_\l(\eta)-r^l_\l(\eta)+\frac{1}{2}\bigl(\Omega^l_\l(\eta)+\Omega^r_\l(\eta)\bigr)
 \quad\mbox{for}
\quad \eta \in \l^*_{op}.
\label{vec}
\ee
Here and further on the elements from $\l \tp \l$ are considered as linear operators from $\l^*$ to $\l$ via
pairing with the first tensor component.

From now on we assume the Lie bialgebra $(\l,\l^*)$ to be factorizable.
This means that the Casimir element $\Omega_\l$ defines an isomorphism $\l^*\to \l$ of $\l$-modules.
Let $\{\xi_i\}\subset \l$ be an orthonormal base with respect to the form $\Omega^{-1}_\l$.
Then the following proposition holds.
\begin{propn}
The quasi-invariant function $r\colon L\to \wedge^2\g$ is a dynamical r-matrix if and only if
the invariant function $r':=r+r_\l$ satisfies the equation
\be
\label{cDYBE1}
\sum_{i} \mathrm{Alt}\bigl(\xi_i\tp \nabla'_{\xi_i}
r'(\la)\bigr)+\mathrm{CYB}\bigl(r'(\la)\bigr)
 &=& Z-Z_l,
\ee
where
$Z_\l:=\CYB(r_\l)$ and
the vector field $\nabla'_\xi$ is defined to be
$\nabla'_\xi:=\frac{1}{2}(\xi^l+\xi^r)$ for $\xi\in \l$.
\end{propn}
\begin{proof}
The elements $r'$ and $r_\l$ satisfy the equality
\be
\CYB(r'-r_\l)&=&\CYB(r')+\CYB(r_\l)-\Alt\bigl([r^{12}_\l+r^{13}_\l,r'_{23}]\bigr)
\nn\\
&=&
\CYB(r')+\CYB(r_\l)+\Alt\bigl(\xi_i\tp r_\l(\eta^i)\tr r'\bigr).
\label{eqCYBsum}
\ee
Here we used invariance of $r'$ with respect to the adjoint action of $\l$.
Now recall that the vector field $r_\l(\eta^i)\tr$ on $L$ equals $r_\l(\eta^i)^l-r_\l(\eta^i)^r$.
It remains to substitute (\ref{eqCYBsum}) into (\ref{cDYBE}) and take into account the specific form of the
$\l^*_{op}$-vector fields on $L$ given by
(\ref{vec}).
\end{proof}
\begin{remark}
Equation (\ref{cDYBE1}) is more general than (\ref{cDYBE}) since it makes sense
for any self-dual Lie algebra $\l$, not necessarily factorizable quasitriangular.
It was derived in \cite{BDF} for the case $\l=\g$.
A solution to this equation for $\l=\g$ (in fact, unique, up to a gauge equivalence)
was found in \cite{FM}.
\end{remark}
\section{Trigonometric dynamical r-matrix over abelian base}
Let $\g$ be a complex  semisimple Lie algebra and $\h$ a Cartan subalgebra in $\g$.
Denote by $\Delta(\g)$ the root system of $\g$ and by the $\Delta_+(\g)$ the subset of positive roots.
Let $\{e_\al\}$, $\al\in \Delta(\g)$ be a base of root vectors from the root subspaces $\g_\al$
normalized to $(e_\al,e_{-\al})=1$ for all $\al\in \Delta_+(\g)$ with respect to
a non-degenerate invariant inner product on $\g$.
Let $\Omega_\g$  denote the split-Casimir of this form,
\be
\Omega_\g:=\sum_{i} x_i\tp x_i+ \sum_{\al\in \Delta(\g)} e_\al\tp e_{-\al},
\nn
\ee
where $\{x_i\}$ is an orthonormal base of $\h$.
Recall from  \cite{EV}, Theorem 3.1, that the meromorphic function
$\h\ni\la\mapsto \rho(\g,\eps,\la)\in \wedge^2 \g$ defined by
\be
\label{Cartan}
\rho(\g,\eps,\la):=\sum_{\al\in \Delta(\g)}
\frac{\eps}{2}\cth\bigl(\frac{\eps}{2}(\al,\la)\bigr)e_\al\tp e_{-\al}
\ee
is a solution to the DYBE with
$Z=Z^\eps_{\Omega_\g}:=\frac{\eps^2}{4}[\Omega^{12}_\g,\Omega^{23}_\g]$.

We may assume $\g$ to be reductive and extend this solution to a meromorphic
function on $\h$ using projection along the center of $\g$.

\section{Base reduction in  trigonometric case}
It is shown in \cite{EV} that the dynamical Yang-Baxter r-matrix admits a reduction of
base when $\l$ is a reductive Lie algebra with the trivial Lie bialgebra structure (abelian $\l^*$).
We will show that an analogous statement also holds
if $\l$ is a factorizable Lie bialgebra.

Let $\l$ be a reductive subalgebra of a complex  Lie algebra $\g$
and let $\h_\l$ be its Cartan subalgebra.
Suppose that $\l$ is equipped
with a factorizable quasitriangular Lie bialgebra structure. Denote by $r_\l\subset \wedge^2\l$ its
classical r-matrix and by $\Omega_\l\in (S^2\l)^\l$ the corresponding Casimir element.
Put the $\l$-base manifold  $L$ to be $\Nc_\l(0)\subset \l$ equipped with the STS Poisson bracket.
Let $\h_\l$ denote the Cartan subalgebra of $\l$.
The open domain $\Nc_{\h_\l}(0):=\Nc_{\l}(0)\cap \h_\l\subset \h_\l$ is a base manifold
for $\h_\l$.
\begin{thm}
\label{reduction}
A quasi-invariant function $r\colon \Nc_{\l}(0)\to \wedge^2\g$  is a dynamical r-matrix over $\Nc_{\l}(0)$ if and only if
the function $\tilde r(\la):=r|_{\h_\l}(\la)+\rho(\la)+r_\l$, where $\rho(\la)=\rho(\g,1,\la)$,
$\la\in  \Nc_{\h_\l}(0)$, is a dynamical r-matrix over the abelian base  $\h_\l$.
\end{thm}
\begin{proof}
Equation (\ref{cDYBE}) on the quasi-invariant function $r$ is equivalent
to equation (\ref{cDYBE1}) on the invariant function $r':=r+r_\l$.
We assume in (\ref{cDYBE1})  that the orthonormal base $\{\xi_k\}$ is compatible with the root decomposition.
Equation (\ref{cDYBE1}) gives rise to the following equation for the  restriction of $r'$ to $\h_\l$:
\newcommand{\rk}{\mathrm{rk}}
\be
\sum_{i=1}^{\rk\> \l} \Alt \Bigl(x_i\tp  \nabla'_{x_i}r'(\la)\Bigr)+\sum_{\al\in \Delta(\l)}
\Alt \Bigl(e_\al\tp \nabla'_{e_{-\al}} r'(\la)\Bigr) + \CYB\bigl(r'(\la)\bigr)=Z_\g-\CYB(r_\l)
\label{eqRestr}
\ee
where $\la\in \h_\l$.
For all $\xi\in \l$ the vector fields $\nabla'_\zeta$,  on $\Nc_{\l}(0)$ read
$$
\nabla'_\xi f(\la)=\partial f(\la)(\frac{1}{2}\ad \la)\cth\bigl(\frac{1}{2}\ad\la\bigr)\xi
,
$$
where $\partial f$ denotes the differential of a function $f$.
On the other hand, $\xi \tr f(\la)=\partial f(\la)(\ad \la)\xi$.
This implies, by $\l$-invariance of $r'$,
$$
e_\al\tp \nabla'_{e_{-\al}} r'(\la)=\cth\bigl(\frac{1}{2}(\al,\la)\bigr)[e_{-\al}\tp 1+ 1\tp e_{-\al},r'(\la)],
\quad\al\in \Delta(\l), \quad \la\in \h_\l.
$$
Observe that $\nabla'_{x_i}=\partial_i:=\frac{\partial}{\partial x_i}$ at $\la\in \h_\l$.
So  equation (\ref{eqRestr}) can be rewritten as
\be
\sum_{i=1}^{\rk\> \l} \Alt \Bigl(x_i\tp  \partial_i r'(\la)\Bigr)+
\Alt \bigl([\rho_{12}+\rho_{13},r'_{23}]\bigr) + \CYB\bigl(r'(\la)\bigr)=Z_\g-\CYB(r_\l)
\nn
\ee
or, upon the substitution $r'=\tilde r-\rho$, as
\be
\sum_{i=1}^{\rk\> \l} \Alt \Bigl(x_i\tp  \partial_i \tilde r(\la)\Bigr)&-&
\sum_{i=1}^{\rk\> \l} \Alt \Bigl(x_i\tp  \partial_i \rho(\la)\Bigr)+
\nn\\
&+& \Alt \bigr([\rho_{12}+\rho_{13},\tilde r_{23}-\rho_{23}]\bigl) + \CYB(\tilde r-\rho)=Z_\g-\CYB(r_\l).
\label{eqRestr1}
\ee
Since $\rho$ is skew, we find $\CYB(\tilde r-\rho)$ to be equal to
\be
&\CYB(\tilde r)+\CYB(\rho)-
([\tilde r_{12},\rho_{13}]+[\tilde r_{12},\rho_{23}]+[\tilde r_{13},\rho_{23}]+
 [\rho_{12},\tilde r_{13}]+[\rho_{12},\tilde r_{23}]+[\rho_{13},\tilde r_{23}])
\nn \\
&
=\CYB(\tilde r)+\CYB(\rho)-\Alt\bigl([\rho_{12}+\rho_{13},\tilde r_{23}]\bigr).
\nn
\ee
Also, it easy to see that $\Alt\bigl([\rho_{12}+\rho_{13},\rho_{23}]\bigr)=2\CYB(\rho)$.
Taking this into account, we rewrite (\ref{eqRestr1})     as
\be
\sum_{i=1}^{\rk\> \l} \Alt \Bigl(x_i\tp  \partial_i \tilde r(\la)\Bigr)
+\CYB(\tilde r)-\sum_{i=1}^{\rk\> \l} \Alt \Bigl(x_i\tp  \partial_i \rho(\la)\Bigr)-\CYB\bigr(\rho(\la)\bigl)=Z_\g-\CYB(r_\l).
\nn
\ee
We have $\CYB(r_\l)=-\CYB(\Omega_\l)=\frac{1}{4}[\Omega^{12}_\l,\Omega_\l^{23}]=Z^{\eps=1}_\l$.
Since $\rho$ solves the CDYBE for $\epsilon=1$, this equation reduces to the
CDYBE for the function
$\tilde r(\la)=r|_{\h_\l}(\la)+\rho(\la)+r_\l$, were $\la$ runs over $\Nc_{\h_\l}(0)$.

Conversely, suppose the function $\tilde r(\la)$ satisfies the CDYBE
over base $\Nc_{\h_\l}(0)$. Then the function $r':=\tilde r|_{\Nc_{\h_\l}(0)}-\rho$ satisfies equation
(\ref{cDYBE1}). Observe that equation (\ref{cDYBE1}) is
$\l$-invariant and, by assumption, $r'$ is the restriction of an invariant function defined on $\Nc_{\l}(0)$.
Therefore the function $r'$ solves (\ref{cDYBE1})
when restricted to semisimple elements from $\Nc_{\l}(0)$ and hence everywhere in $\Nc_{\l}(0)$
 since semisimple elements are dense in $\l$.
\end{proof}

Theorem \ref{reduction} allows to construct a trigonometric dynamical r-matrix
over a Levi subalgebra in a complex simple Lie algebra.
\begin{corollary}
\label{q-gr}
Let $\l$ be a Levi subalgebra in a complex simple Lie algebra $\g$ and let $\l^\bot$ be the
orthogonal complement of $\l$ in $\g$ with respect to the Killing form.
The meromorphic function
$r\colon \Nc_{\l}(0)\ni \la\mapsto f_\l(\ad \la)\oplus f_{\l^\bot}(\ad \la)\in \End(\l)\oplus \End(\l^\bot)\subset \End(\g)$ with
\be
f_{\l}(\la) &=&\frac{1}{2}\cth\bigl(\frac{1}{2}\ad \la \bigr)-
\frac{\eps}{2}\cth\bigl(\frac{\eps}{2}\ad \la \bigr)-r_\l,
\nn
\\
f_{\l^\bot}(\la)&=& -\frac{\eps}{2}\cth\bigl(\frac{\eps}{2}\ad \la \bigr).
\nn
\ee
is a PL dynamical r-matrix for $Z^\eps_\g=\frac{\eps^2}{4}[\Omega^{12}_\g,\Omega^{23}_\g]$.
\end{corollary}
\begin{proof}
Since $\g$ is simple, the restriction of the matrix $\ad \la\in \End(\g)$ to the invariant subspace
$\l^\bot$ is invertible for $\la$ belonging to a dense open subset in $\Nc_{\l}(0)$.
Hence  $\la\mapsto r(\la)$ is a correctly defined  meromorphic function on $\Nc_{\l}(0)$.
The function $r+r_\l$ is $\l$-invariant and $r|_{\h_\l}+r_\l+\rho(\l,1,\la)=\rho(\g,\eps,\la)$
for $\la\in \h_\l$ .
\end{proof}

\section{Generalized Feh\'er-Marshall dynamical r-matrices}
Let $\g$ be a finite dimensional complex Lie algebra and $B\colon \g\to \g$ an automorphism
of order $n$. Then $\g=\oplus_{j\in \Z/n\Z}\g_j$, where $\g_j:=\ker (B-e^{2ipj/n})$.
The Lie algebra $\g_0$ acts on $\g_j$ for all $j$.

Let $\g$  be equipped with an $\ad$- and $B$-invariant form; denote by $\Omega_\g\in (S^2\g)^\g$ the
corresponding Casimir element. According to
\cite{ES}, there exists a dynamical r-matrix over the base $\l^*$ for
the trivial Lie bialgebra $\l$ (abelian $\l^*$). Under the identification $\g\tp \g\simeq \g^*\tp \g\simeq \End(\g)$, it is given
by the invariant function $\rho\colon \mathcal{N}(\l)\to \End(\g)$, $\rho(A)|_{\g_j}=f_j(\ad A)$, with
\be
f_0(s)&=&\frac{1}{s}-\frac{1}{2}\cth\bigl(\frac{1}{2}s\bigr),
\label{r-tr1}\\
f_j(s)&=&-\frac{1}{2}\cth\bigl(\frac{1}{2}s+\frac{i\pi j}{n}\bigr),\quad j\not = 0.
\label{r-tr2}
\ee
This solution corresponds to $Z_\g=\frac{1}{4}[\Omega^{12}_\g,\Omega^{23}_\g]$.
For $B=\id$, formulas (\ref{r-tr1}-\ref{r-tr2}) give the rational-trigonometric dynamical r-matrix
of Alekseev-Meinrenken, \cite{AM}.

The trigonometric analog of this construction is as follows.
Consider the analytical function $r'_{trig}\colon \mathcal{N}(\l)\to \End(\g)$, $r'_{trig}(A)|_{\g_j}=f^{trig}_j(\ad A)$, with
\be
\label{tr1}
f^{trig}_0(s)&=&\frac{1}{2}\cth\bigl(\frac{1}{2} s\bigr)-\frac{\eps}{2}\cth\bigl(\frac{\eps}{2} s\bigr),
\\
f^{trig}_j(s)&=&-\frac{\eps}{2}\cth\bigl(\frac{\eps}{2}s+\frac{i\pi j}{n}\bigr),\quad j\not = 0.
\label{tr2}
\ee
Let $\Omega_\l$ denote the Casimir element of the restriction of the invariant form
to $\l$.
\begin{thm}
The equivariant function $r'_{trig}$ satisfies equation (\ref{cDYBE1})
for $Z^\eps_\g=\frac{\eps^2}{4}[\Omega^{12}_\g,\Omega^{23}_\g]$
and $Z_\l=\frac{1}{4}[\Omega^{12}_\l,\Omega^{23}_\l]$.
\end{thm}
\begin{proof}
The proof is an appropriate modification of the proof of Theorem A.1 of \cite{ES}.
At the first step one considers $\g=\l\oplus\ldots \oplus \l $, the direct sum of $n$-copies of a reductive
Lie algebra $\l$ and $B$ the cyclic permutation of these copies. One proves that the
function $r'_{trig}|_{\h_l}+\rho$, where $\rho(\la):=\rho(\l,1,\la)$ from (\ref{Cartan})
is a dynamical r-matrix over the $\h_l$
for $Z=Z^\eps_\g$. At $\eps=1$, this function coincides with the corresponding function
from \cite{ES}, see the proof of Proposition A.1; for $\eps\not =1$ cf. Remark \ref{dilation}.
 By Theorem \ref{reduction}, $r'$ solves equation (\ref{cDYBE1}) for these specific $\g$ and $B$.
The case of general $\g$ is derived from this one similarly to \cite{ES}.
Namely, define a map $W\colon \Nc_{\l}(0)\to \wedge^3\g$ setting
$$
W(A)=\Alt\bigl(\nabla'r'_{trig}(A)\bigr)+\CYB\bigl(r'_{trig}(A)\bigr)-Z^\eps_\g+\frac{1}{4}[\Omega^{12}_\l,\Omega^{13}_\l].
$$
This is an analog of the function $W$ from Appendix of \cite{ES}.
Then one defines
$K_{ij}(A,X,Y):=\bigl(\id\tp X\tp Y,W(A)\bigr)\in \g_{i+j}$ for all $i,j\in \Z/n\Z$, $A \in \l$, $X\in \g_i$ and $Y_j\in \g_j$,
which are expressed by universal Lie series in $A$, $X$, and $Y$ for each pair $i,j$.
For $K_{ij}(A,X,Y)$ one can prove the same assertions as Proposition A.2 and Lemma A.2 of \cite{ES}
and find that they vanish for all  $A\in \l$, $X\in \g_i$, and $Y_j\in \g_j$.
Thus $W(A)=0$ and this proves the theorem.
\end{proof}
\begin{corollary}
Suppose $\l=\ker B$ has a quasitriangular structure $(r_\l,\Omega_\l)$.
Then $r_{trig}=r'_{trig}-r_\l$ is a PL dynamical r-matrix.
\end{corollary}
For $B=\id$, formulas  (\ref{tr1}-\ref{tr2}) give the trigonometric PL dynamical r-matrix of
Feh\'er-Marshal, \cite{FM}.

\end{document}